\documentclass[11pt]{article}
\PassOptionsToPackage{obeyspaces}{url}
\usepackage{amsfonts, amsmath, amssymb}

\usepackage[hidelinks]{hyperref}
\usepackage{bookmark}
\bookmarksetup{open,numbered,depth=5} 

\usepackage{amsthm}

\usepackage{breakurl}
\usepackage{tikz}
\usetikzlibrary{calc,patterns,tikzmark}

\usepackage{etoolbox} 
\patchcmd{\thebibliography}{\leftmargin\labelwidth}{\leftmargin\labelwidth\addtolength\itemsep{-0.3\baselineskip}}{}{}

\usepackage{mathtools}

\usepackage[nameinlink]{cleveref}

\usepackage{tabularx}

\newtheorem{theorem}{Theorem}
\newtheorem{lemma}[theorem]{Lemma}
\newtheorem{corollary}[theorem]{Corollary}
\newtheorem{proposition}[theorem]{Proposition}

\crefname{conjecture}{Conjecture}{Conjectures}

\newcommand*\samethanks[1][\value{footnote}]{\footnotemark[#1]} 
\author{Boris Bukh\thanks{Department of Mathematical Sciences, Carnegie Mellon University, Pittsburgh, PA 15213, USA\@. Supported in part by U.S.\ taxpayers through NSF grant DMS-2154063.}\,\,\textsuperscript{,}\thanks{Supported in part by a Simons Foundation Fellowship.} \and Aleksandre Saatashvili\samethanks[1]}

\title{Maximal sets of a given diameter in Hamming cubes}

\date{}

\newcommand*{\eqdef}{\stackrel{\mbox{\normalfont\tiny def}}{=}}   

\newcommand*{\N}{\mathbb{N}}                                     
\newcommand*{\E}{\mathop{\mathbb{E}}}                            
\newcommand*{\R}{\mathbb{R}}                                     
\DeclarePairedDelimiter\abs{\lvert}{\rvert}                     
\DeclarePairedDelimiter\norm{\lVert}{\rVert}                     

\DeclareMathOperator{\dist}{dist}                               
\DeclareMathOperator{\diam}{diam}                               
\DeclareMathOperator{\wt}{wt}                                   
\DeclareMathOperator{\Sun}{Sun}                               
\DeclareMathOperator{\Ball}{Ball}                               
\newcommand*{\Ninfty}{[n]^\infty}
\newcommand*{\Tinfty}{([n]\cup\{\star\})^\infty}
\newcommand*{\twoinfty}{[2]^\infty}
\newcommand*{\Pbig}{P_{\text{big}}}
\newcommand*{\Psmall}{P_{\text{small}}}

\begin{document}
\maketitle

\begin{abstract}
  A subset of the Hamming cube over $n$-letter alphabet is said to be \emph{$d$-maximal} if its diameter is $d$, and adding any point increases the diameter.
  Our main result shows that each $d$-maximal set is either of size at most $(n+o(n))^d$ or contains a non-trivial Hamming ball. The bound of $(n+o(n))^d$
  is asymptotically tight. Additionally, we give a non-trivial lower bound on the size of any $d$-maximal set and show that the number of essentially different
  $d$-maximal sets is finite.
\end{abstract}

\section*{Introduction}
\paragraph{Motivation.}
Consider the Boolean cube $\{0,1\}^N$ endowed with the usual Hamming distance.
A set $S$ is \emph{$d$-maximal} if its diameter is $d$, and it is an inclusion-maximal
set with such a property.

We can lift $S$ to a subset of $\{0,1\}^{N'}$, for $N'>N$, by appending the same suffix to all the words in~$S$.
The resulting set still has diameter $d$, and so is contained in a (possibly larger)
$d$-maximal set. This leads to the natural question, asked by Briggs, Feng and Wells \cite{briggs_cox_feng}, of whether there exist arbitrarily large
$d$-maximal sets that are stable under lifting. Because the $d$-maximal subsets of $\{0,1\}^N$ that are stable under lifting correspond to \emph{finite} $d$-maximal subsets of $\{0,1\}^{\infty}$,
these are the objects that we study in this paper.

\paragraph{Upper bounds on the size of a finite \texorpdfstring{$d$-maximal}{d-maximal} set.}
Because our arguments work not only for $\{0,1\}$, but for an arbitrary alphabet,
we state our results in this generality. As the nature of the alphabet is not important to us,
we use $[n]\eqdef \{1,2,\dotsc\}$ as our standard alphabet.

We denote by $\dist(a,b)$ the \emph{Hamming distance} between $a,b\in \Ninfty$.
For a set $S\subset \Ninfty$,
its \emph{diameter} is
\[
  \diam(S)\eqdef \max_{a,b\in S}\dist(a,b).
\]
A set $S\subset \Ninfty$ is \emph{$d$-maximal} if
its diameter is $d$, and it is largest (under inclusion) set of diameter~$d$.

There are lots of infinite $d$-maximal sets.
For example, for even $d$, any Hamming ball of radius $d/2$ is $d$-maximal. More generally, using
Zorn's lemma one can show that any set
of diameter at most $d$ can be extended to a $d$-maximal set; in particular, any
infinite set of diameter at most $d$ is contained in an infinite $d$-maximal set.

We prove that, not only finite $d$-maximal sets have bounded size, but the same
holds for $d$-maximal sets containing no Hamming ball of radius $1$
(or \emph{$1$-balls} for brevity).

 \begin{theorem}\label{thm:main}
     Suppose that $S\subset \Ninfty$ is a $d$-maximal set that contains no $1$-ball. Then
     $$|S|\leq d^2(n+8n^{2/3})^d.$$
\end{theorem}
The theorem implies the same result about subsets of $[n]^M$ for finite $M$. Indeed, suppose $S\subset [n]^M$ is a $d$-maximal set of size
more than $d^2(n+8n^{2/3})^d$, for some finite $M$. Embed $S$ into $\Ninfty$ and extend it to a $d$-maximal set $\hat{S}\subset \Ninfty$.
By \Cref{thm:main} the set $\hat{S}$ contains a $1$-ball around some point $a\in \hat{S}$. Let $a_{\leq M}\eqdef (a_1,\dotsc,a_M)$.
Maximality of $S$ then implies that the $1$-ball around $a_{\leq}$ is contained in~$S$. 

The bound in \Cref{thm:main} is nearly sharp, for the
the $d$-dimensional cube $[n]^d \times\nobreak \{0\}$ has
$n^d$ elements and is $d$-maximal.
We suspect that the cube is the largest finite $d$-maximal
set for sufficiently large $n$. However, for $n=2$, larger sets exist.
\begin{theorem}\label{thm:binaryconstr}
There exist $d$-maximal sets of size more than $\binom{\lfloor 3d/2\rfloor}{d}\geq \bigl(2.59+o(1)\bigr)^d$ in $\twoinfty$, for all~$d$.
\end{theorem}
The best upper bound we have been able prove for $n=2$ is the following.
\begin{theorem}\label{thm:binary}
    Every $d$-maximal set $S\subset \twoinfty$ not containing
    a $1$-ball is of size $$|S|\leq  (4-10^{-10})^d,$$
    for all large $d$.
\end{theorem}

\paragraph{More general upper bounds and the number of distinct \texorpdfstring{$d$-maximal}{d-maximal} sets.}
One can go beyond \Cref{thm:main} and ask about the structure infinite $d$-maximal sets.
For example, suppose that $S$ is an infinite $d$-maximal set, and $S'$ is obtained from $S$ by removing
points of an arbitrarily chosen $1$-ball from~$S$; must the leftover set $S'$ be of bounded size if it is finite? It turns out that
not only the answer to this question is yes, but much more is true. 

Denote by $\Ball(a,\ell)$ the $\ell$-ball around point $a$. Let
\[
  S^{(\ell)}\eqdef \{a\in S : \Ball(a,\ell)\subset S,\ \nexists b\, \Ball(a,\ell)\subset \Ball(b,\ell+1)\subset S\}.
\]
Note that if $S$ is $1$-ball-free, then $S^{(0)}=S$. So, the following is a generalization of \Cref{thm:main}.
\begin{theorem}\label{thm:mainL}
  Suppose that $S\subset \Ninfty$ is a $d$-maximal set that contains no $(L+1)$-ball. Then
  \[
    \abs{S^{(\ell)}}\leq Cd^{2L+2}(n+8n^{2/3})^d.
  \]
  for some constant $C=C(L,n)$ that depends only on $L$ and~$n$.
\end{theorem}
As a corollary of this bound, we can show that there are only finitely many distinct $d$-maximal sets in $\Ninfty$,
up to permutation of the alphabet (separately for each position) and permutation of the coordinates.
Formally, a map $\sigma\colon \Ninfty\to \Ninfty$ is an \emph{isomorphism} if there are permutations $\sigma_1,\sigma_2,\dotsc$ of $[n]$
and a permutation $\bar{\sigma}$ of $\N$ such that $\sigma(a)_i=\sigma_i(a_{\bar{\sigma}(i)})$.
\begin{corollary}
There are only finitely many isomorphism types of $d$-maximal sets in $\Ninfty$.
\end{corollary}
\begin{proof}
  Let $S$ be a $d$-maximal set. Set $N\eqdef \sum_{\ell\leq d} \abs{S^{(\ell)}}$.
  By applying a suitable automorphism we may ensure that the all-$1$ word is in $S$, and so
  there are at most $d$ non-$1$ symbols in any $a\in S$. After permuting the coordinates,
  we may also assume that $a_i=1$ for all $i>dN$ and all $a\in \bigcup_{\ell\leq d} S^{(\ell)}$.
  Since $S=\bigcup_{\ell\leq d}\bigcup_{a\in S^{(\ell)}} \Ball(a,\ell)$, and $N\leq G(d)$ for some function $G$ by
  \Cref{thm:mainL}, the result follows.
\end{proof}

\paragraph{Lower bounds on the size of a finite \texorpdfstring{$d$-maximal}{d-maximal} set.} We saw that the size of the largest $d$\nobreakdash-maximal set grows exponentially in $d$.
In contrast, \emph{smallest} $d$-maximal set are much smaller. Briggs, Feng and Wells \cite{briggs_cox_feng} constructed a $d$-maximal set of size $2d$ in $\{0,1\}^\infty$, assuming
the existence of a Hadamard matrix of order $2d$. It is conjectured that Hadamard matrices of order $4k$ exist for all positive integers $k$, and the conjecture
is known for infinitely many values of $k$.

We give a polynomial lower bound on the size of a $d$-maximal set.
\begin{theorem}\label{thm:lower}
Every $d$-maximal set has at least $\tfrac{1}{2}(d/\log 2d)^{2/3}$ elements.
\end{theorem}

We suspect that the construction in \cite{briggs_cox_feng} is optimal, and the smallest $d$-maximal sets
have size at least $2d$. This would be the best possible, not only for $n=2$, but for larger alphabets as well.
This is a consequence of the following result, which shows that most constructions, including the one from~\cite{briggs_cox_feng},
extend to higher alphabets.
\begin{proposition}
  Assume that $S\subset [n]^r$ is a $d$-maximal set such that for every $a\in S$ and every coordinate $i\in [r]$
  there is a $b\in S$ such that $a_i=b_i$ and $\dist(a,b)=d$. Then
  $S$ is also $d$-maximal as a subset of $[n+1]^r$.
\end{proposition}
\begin{proof}
  Suppose that $S$ is not maximal in $[n+1]^r$, and $a'\notin S$ satisfies $\dist(a',S)\leq d$.
  Let $a\in [n]^r$ be obtained from $a'$ by replacing all occurences of $n+1$ by~$n$.
  Then $\dist(a,c)\leq \dist(a',c)\leq d$ for all $c\in S$, and so $a\in S$ by maximality of~$S$.
  Let $b$ be an element satisfying $a_i=b_i$ and $\dist(a,b)=d$. Then $\dist(a',b)>\dist(a,b)=d$,
  contradicting $\dist(a',S)\leq d$.
\end{proof}

\paragraph{Acknowledgement.} We are thankful to Chris Wells for valuable discussions.
\section{A weak upper bound on the sizes of sets \texorpdfstring{$S^{(\ell)}$}{S\textasciicircum(ell)}}\label{sec:basic}
In this section we prove a version of \Cref{thm:mainL} with a weaker bound. We then later use this bound as an ingredient in the 
proofs of stronger results: the large-alphabet bound of \Cref{thm:mainL} 
in the next section, and the binary alphabet bound of \Cref{thm:binary} in \Cref{sec:binary}.\medskip

We call $t\in \Tinfty$ a \emph{template}. We regard elements of $[n]$ as fixed and $\star$ as a wildcard.
We say that the word $a\in \Ninfty$ \emph{fits} the template $t$ if $a_i=t_i$ whenever $t_i\neq \star$.
For a set $X\subset \Ninfty$ and a template $t$ let
\[
  X[t]\eqdef \{ a\in X : a\text{ fits }t\}.
\]

Let $S$ be a $d$-maximal set not containing any $(L+1)$-ball. We may assume without loss of generality that $L\leq d/2$.
Let $w\in S$ be a  word in~$S$.
We say that a template $t$ is \emph{regular with respect to $w$} if $t_i\neq w_i$ for all~$i$.
When the word $w$ is fixed, we say simply that $t$ is \emph{regular}.
The \emph{weight} of $t$ 
is $\wt(t)\eqdef \abs{\{i : t_i\neq \star\} }$.

To prove the simpler bound, let $S'\subset S$ be a set whose size we wish to bound, e.g., in \Cref{thm:mainL} it could be $S^{(\ell)}$.
We then fix $w\in S$ arbitrarily, and build a sequence of regular templates $t^{(0)},t^{(1)},\dotsc$ of weight $\wt(t^{(m)})=m$
such that each $S'[t^{(m+1)}]$ is not much smaller than $S'[t^{(m)}]$. Because a regular
template of weight $m$ satisfies $\dist(w,S'[t^{(m)}])\geq m$, the
sets $S'[t^{(m)}]$ must eventually be rather small. Taken together, these two properties of $S'[t^{(m)}]$ will
imply an upper bound on~$\abs{S'}$,

For ease of bookkeeping, we shall work with subsets of $S_k\eqdef \{a\in S : \dist(w,a)=k\}$ rather than with $S$.
The following lemma is the main inductive step.

\begin{proposition}\label{prop-induct} 
  Suppose that $t^{(m)}$ is a regular template of weight $m$ and $S'\subset S_k[t^{(m)}]$.
  Then there exists a regular template $t^{(m+1)}$ of weight
  $m+1$ satisfying
  \[
    \abs{S'[t^{(m+1)}]}\geq \frac{1}{n}\cdot \frac{k-m-L}{d-m/2-L/2} \abs{S'}.
  \]
\end{proposition}
In this lemma, the set $S'$ is allowed to be infinite. If this is the case, then the proposition asserts that $S'[t^{(m+1)}]$ is infinite, provided that $k>m+L$.

We need the following lemma, which is the only place in the proof where we use the assumption that the $d$-maximal
set $S$ contains no $(L+1)$-ball.
\begin{lemma}\label{max-lemma}
For every $y\in \Ninfty$ there exists $r\in S$ such that $\dist(r,y)\geq d-L$.
\end{lemma}
\begin{proof}
  Assume the contrary. Consider any point $z$ satisfying $\dist(z,y)\leq L+1$. By
  the triangle inequality the $d$-ball around $z$ contains $S$. By maximality of $S$, that entails $z\in S$.
  Since $z$ is arbitrary, this means that $S$ contains the $(L+1)$-ball around $y$, contradicting
  the assumption on~$S$.
\end{proof}

\begin{proof}[Proof of \Cref{prop-induct}]
We may assume that $m\leq k\leq d$, for otherwise the set $S_k[t^{(m)}]$ is empty, rendering the result trivial.

Let $y$ be the word defined by
\[
  y_i\eqdef
  \begin{cases}
    t_i&\text{if }t_i\neq \star,\\
    w_i&\text{if }t_i=\star.
  \end{cases}
\]
By \Cref{max-lemma} there is $r\in S$ such that $\dist(r,y)\geq d-L$. Let
\begin{align*}
  \alpha &\eqdef \{ i : r_i\neq y_i,\  t_i\neq \star\},\\
  \beta &\eqdef \{ i : r_i\neq y_i,\  t_i= \star\}.\\
  \intertext{For $z\in S'$, define}
  \gamma(z) &\eqdef \{ i\in \beta : z_i=r_i\},\\
  \delta(z) &\eqdef \{ i\in \beta : z_i\notin \{r_i,w_i\}\}.
\end{align*}
These definitions are illustrated in the following figure (where the values of $w,t,y,r,z$ are chosen as samples).
\begin{align*}
    w\ &1111\overbrace{1111}^{\alpha}\overbrace{1111\,1111\,1111}^{\beta}\,1111\\
    t\ &2222\,2222\,{\star}{\star}{\star}{\star}\,{\star}{\star}{\star}{\star}\,{\star}{\star}{\star}{\star}\,{\star}{\star}{\star}{\star}\\  
    y\ &2222\,2222\,1111\,1111\,1111\,1111\\
    r\ &2222\,1133\,2222\,2222\,2222\,1111\\
    z\ &\tikzmarknode{A}{2}222\,222\tikzmarknode{B}{2}\,1111\underbrace{2222}_{\gamma(z)}\underbrace{3333}_{\delta(z)}\,????
\begin{tikzpicture}[overlay, remember picture]
\draw[<->] ($(A.south west)+(0,-.3)$) -- ($(B.south east)+(0,-.3)$) node [midway,below=-2pt]{$\scriptstyle m$};
\end{tikzpicture}
\end{align*}

Note that
\begin{align}
  d-L\leq \dist(r,y)&=\abs{\alpha}+\abs{\beta},\label{eq-ry}\\
  d\geq \dist(r,w)&\geq (m-\abs{\alpha})+\abs{\beta}.\label{eq-rw}\\
\intertext{Furthemore, since $\wt(t^{(m)})=m$ and $z\in S_k$, there are exactly
  $k-m-\gamma(z)-\delta(z)$
  positions $i\notin\alpha\cup\beta$ satisfying $y_i\neq z_i$. So, breaking the contribution to $\dist(r,z)$ from
  coordinates in $\alpha$, in $\beta$ and the rest, we obtain}
  \label{eq-rz}
  d\geq \dist(r,z)&=\abs{\alpha}+\bigl(\abs{\beta}-\abs{\gamma(z)}\bigr)+\bigl(k-m-\abs{\gamma(z)}-\abs{\delta(z)}\bigr).
\end{align}
Rearranging this gives
\begin{equation}\label{gamma-delta}
  2\abs{\gamma(z)}+\abs{\delta(z)}\geq k-m+(\abs{\alpha}+\abs{\beta}-d).
\end{equation}
If the set $S'$ is infinite, we may conclude the proof at this point by noting
that the right side of \eqref{gamma-delta} is positive and both $\gamma(z)$ and $\delta(z)$ are subsets of
a finite set~$\beta$.

So, assuming that $S'$ is finite, pick the word $z$ uniformly at random from $S'$, and define
\[
  p_i(\ell)\eqdef \Pr [z_i=\ell].
\]
We have
\begin{align*}
  \sum_{i\colon i\in \beta}p_i(r_i)&=\E \abs{\gamma(z)},\\
  \sum_{\substack{i\colon i\in \beta\\\ell\notin\{ w_i,r_i\}} }p_i(\ell)&=\E \abs{\delta(z)},\\
\intertext{which implies, by means of \eqref{gamma-delta}, that}
  \sum_{i\colon i\in \beta}\left(2 p_i(r_i)+\sum_{\ell\notin \{w_i,r_i\}} p_i(\ell)\right)&\geq k-m+(\abs{\alpha}+\abs{\beta}-d).
\end{align*}
Hence,
\begin{align}\label{max-prob}
  \max_{\substack{i\colon i\in \beta\\\ell\neq w_i}} p_i(\ell)&\geq \frac{k-m+(\abs{\alpha}+\abs{\beta}-d)}{n\abs{\beta}}\stackrel{\eqref{eq-rw}}{\geq}\frac{1}{n}\cdot
  \frac{k-m+(\abs{\alpha}+\abs{\beta}-d)}{d-m/2+(\abs{\alpha}+\abs{\beta}-d)/2}\\\notag
  &\stackrel{\eqref{eq-ry}}{\geq} \frac{1}{n}\cdot \frac{k-m-L}{d-m/2-L/2},
\end{align}
where in the last inequalty we used that $L\leq d/2$ and $m\leq d$ together imply that $d-m/2-L/2>0$.

Pick any pair $(i,\ell)$ attaining maximum in \eqref{max-prob}, and let $t^{(m+1)}$ be obtained from $t^{(m)}$ by replacing its $i$'th letter with $\ell$.
Since $\wt(t^{(m+1)})=\wt(t^{(m)})+1$ and $\abs{S'[t^{(m+1)}]}=p_i(\ell)\abs{S'}$, the lemma follows.
\end{proof}

The preceding lemma on its own is enough to obtain a bound for $\abs{S_k}$ in the case $L=0$; one applies it inductively
to patterns of ever-increasing weight and uses the fact that $\abs{S_k[t^{(k)}]}\leq 1$ for any pattern of weight~$k$.
To treat the general case, we need a generalization of the latter fact to sets $S^{(\ell)}_k$ for any~$L$.

A \emph{$p$-sunflower} is a family of $p$ sets whose pairwise intersections are identical.
Erd\H{o}s and Rado \cite{erdos_rado} proved that any sufficiently large family of $k$-element sets contains a $p$-sunflower.
The following is the current best bound, due to Alweiss, Lovett, Wu, Zhang \cite{alweiss_lovett_wu_zhang} with the slight improvement by Bell, Chueluecha, Warnke \cite{bell_chieluecha_warnke}.
\begin{lemma}[Sunflower lemma]\label{sunflower}
  For any $p,k$ there is a number $\Sun(p,k)$ such that every
  family of more than $\Sun(p,k)$ many distinct $k$-element sets contains a $p$-sunflower,
  where
  \[
    \Sun(p,k)\leq (Cp\log k)^k\qquad\text{ for }p,k\geq 2
  \]
  and $C>4$ is an absolute constant.
\end{lemma}
More generally, we say that the words $a^1,\dotsc,a^p\in \Ninfty$ 
form a \emph{$p$-sunflower with respect to $w$} if, for each position $i$, either
$a^1_i=\dotsb=a^p_i$ or there is a unique $j$ satisfying $a^j_i\neq w_i$.
We call $\{i : a^1_i=\dotsb=a^p_i\neq w_i\}$ the \emph{stem} of the sunflower.
\begin{lemma}[Generalized sunflower lemma]
  Suppose that $p\geq n$. Then any family $X\subset \Ninfty$ of more than $\Sun(p,2k)$ many distinct words at Hamming distance $k$
  from $w$ contains a $p$-sunflower with respect to $w$.
\end{lemma}
\begin{proof}
  Without loss of generality, $w_i=n$ for all $i$. 
  For a word $a\in X$, define sets
  \begin{align*}
    C_1(a)&\eqdef \{ i\in \N : a_i\neq n\},\\
    C_2(a)&\eqdef \{ (i,\ell)\in \N\times [n-1] : a_i=\ell\},
  \end{align*}
  and let $C(a)\eqdef C_1(a)\sqcup C_2(a)$ be their disjoint union.
  It is clear that the map $a\mapsto C(a)$
  is a bijection and that $\abs{C(a)}=2k$ for all $a\in X$. By the sunflower lemma, there are $p$ words $a^1,\dotsc,a^p\in X$
  such that the sets $C(a^1),\dotsc,C(a^p)$ form a $p$-sunflower. We claim that the words $a^1,\dotsc,a^p$ form a $p$-sunflower with respect to~$w$.
  Indeed, let $i\in \N$ and suppose that neither $a^{j}_i$ nor $a^{j'}_i$ is equal to~$w_i$. Then $i\in C_1(a^{j})\cap C_1(a^{j'})$
  implying that $a^k_i\neq w_i$ for all $k\in [p]$. Since $p>n-1$, there are two words among the $a^1,\dotsc,a^p$
  that have the same $i$'th coordinate. Since $C_2(a^1)=\dotsb=C_2(a^p)$, this in turn implies that all of these words
  have the same $i$'th coordinate.
\end{proof}

Recall that $S^{(\ell)}$ was defined as the set of centers of $\ell$-balls in $S$ that are not contained in
any larger ball. Let $S^{(\ell)}_k\eqdef S^{(\ell)}\cap S_k$.
\begin{lemma}\label{sunflower-bound}
  Let $t^{(m)}$ be a regular template of weight $m$. Then
  \[
    \abs{S^{(\ell)}_k[t^{(m)}]}\leq \Sun(d+n,2k-2m).
  \]
\end{lemma}
Note that, pedantically speaking, $\Sun(p,0)=1$ because there exists only one $0$-element set. In this case, the
lemma correctly asserts that at most one element of $S^{(\ell)}_k$ fits a pattern of weight~$k$.
\begin{proof}[Proof of \Cref{sunflower-bound}]
  Let $p\eqdef d+n$, and assume, on the contrary, that $\abs{S^{(\ell)}_k[t^{(m)}]}> \Sun(p,2k-2m)$.
  By renumbering the coordinates, we may assume without loss of generality that the first $m$ symbols of $t$
  are in $[n]$, and the rest are all $\star$'s. This way, each word $a\in S^{(\ell)}_k[t^{(m)}]$
  begins with the same prefix of length~$m$. Denote by $a_{>m}$ the suffix of $a$ obtained by deleting the first $m$ symbols.
  Since $t^{(m)}$ is a regular pattern, $\dist(a_{>m},w_{>m})=\dist(a,w)-m=k-m$. Therefore, by the generalized sunflower lemma
  there are $a_1,\dotsc,a_p\in S^{(\ell)}_k[t^{(m)}]$ such that $a^1_{>m},\dotsc,a^p_{>m}$ form a $p$-sunflower. Since $a^1,\dotsc,a^p$ all share the same prefix,
  this means that $a^1,\dotsc,a^p$ themselves form a $p$-sunflower.

  By renumbering the coordinates again, we may assume that $[m']$ is the sunflower's stem, for some $m'\geq m$.
  Define a word $b$ by
  \[
    b_i\eqdef
    \begin{cases}
      a^1_i&\text{for }i\leq m',\\
      w_i&\text{for }i>m.
    \end{cases}
  \]
  We claim that $\Ball(b,\ell+k-m')\subset S$. Indeed, let $z\in S$ be arbitrary; we shall show the containment $\Ball(b,\ell+k-m')\subset \Ball(z,d)$.
  Because the sets of positions where $a^j_{>m'}$ differs from $w_{>m'}$ are all disjoint, and because
  $\dist(w,z)\leq d<p$, there is an index $j$ such that the set $\{i: a^j_i=z_i\neq w_i\}$ is contained in $[m']$. Therefore,
  \begin{align*}
    \dist(a^j,z)&=\dist(a^j_{\leq m'},z_{\leq m})+\dist(a^j_{> m'},z_{> m})\\
                &=\dist(b_{\leq m'},z_{\leq m'})+\dist(w_{>m'},z_{>m'})+\dist(w_{>m'},a^j_{>m'})\\
                &=\dist(b_{\leq m'},z_{\leq m})+\dist(b_{>m'},z_{>m'})+k-m'\\
                &=\dist(b,z)+k-m'.
  \end{align*}
  Since $S$ is a diameter-$d$ set, $\Ball(a^j,\ell)\subset \Ball(z,d)$, which by the triangle inequality implies the containment
  $\Ball(b,\ell+k-m')\subset \Ball(z,d)$. Since $z\in S$ is arbitrary and the set $S$ is $d$-maximal, it follows that
  $\Ball(b,\ell+k-m')\subset S$.

  Since $\Ball(a^j,\ell)\subset \Ball(b,\ell+k-m')$ and $a^j\in S^{\ell}$, this contradicts the definition of $S^{\ell}$
  unless $m'=k$. The case $m'=k$ is not possible as well, as in this case the condition $\dist(a^j,w)=k$ means that $a^1=\dotsb=a^p$,
  contradicting $p>1$.
\end{proof}

\begin{corollary}\label{Sk-template-bound}
  Suppose that $t^{(m)}$ is a regular template of weight $m$. 
  Then
  \[
    \abs{S^{(\ell)}_k[t^{(m)}]}\leq \Sun(d+n,2L)(n/2)^{k-m-L}\binom{2d-m-L}{k-m-L}.
  \]
\end{corollary}
\begin{proof}
  By induction on~$k-m-L$. The base case $m=k-L$ follows from \Cref{sunflower-bound}.
  The induction step follows from \Cref{prop-induct} applied to $S'=S_k^{(\ell)}[t^{(m)}]$.
\end{proof}

\begin{corollary}\label{weak-bound}
  Suppose that $S$ is a $d$-maximal set that contains no $1$-ball. Then
  \[
    \abs{S}\leq (2n)^d.
  \]
\end{corollary}
\begin{proof}
  Since $S_k[t^{(0)}]=S_k$, it follows from the above that
  \[
    \abs{S_k}\leq n^k\frac{d}{k}\cdot\frac{d-1/2}{k-1}\cdot \dotsb \frac{d-k/2+1/2}{1}=(n/2)^k\binom{2d}{k}.
  \]
  Hence,
  \[
    \abs{S}=\sum_{k=0}^d \abs{S_k}\leq \sum_{k=0}^d (n/2)^k\binom{2d}{k}\leq (n/2)^d\sum_{k=0}^d \binom{2d}{k}\leq 4^d(n/2)^d.\qedhere
  \]
\end{proof}

\section{Sizes of \texorpdfstring{$d$-maximal}{d-maximal} sets over large alphabet: proof of \texorpdfstring{\Cref{thm:mainL}}{Theorem 4}}
We may assume that $n\geq 100$, for otherwise $n+8n^{2/3}\leq 2n$ and \Cref{thm:main} follows from \Cref{weak-bound}.

Let $a,b\in S$ be any pair of points at distance $d$ apart; these exist by maximality of $S$.
By reordering the coordinates, we may assume that $a$ and $b$ differ in the first $d$ coordinates.
For a word $z\in \Ninfty$, let $z_{\leq d}\in [n]^d$ be the prefix of $z$ of length~$d$.
For $p\in [n]^d$ define
\begin{align*}
  S^{(\ell)}(p)&\eqdef \{ z\in S^{(\ell)} : z_{\leq d}=p\},\\
  m(p)&\eqdef \max\bigl(\dist(a_{\leq d},p),\dist(b_{\leq d},p)\bigr).
\end{align*}

For sake of notational brevity, let $Z\eqdef \Sun(d+n,2L)$. 
\begin{lemma}
  For every $p\in [n]^d$ we have $\abs{S^{(\ell)}(p)}\leq d Z (n/2)^{d-m(p)} \binom{2d-m(p)}{d-m(p)}$.
\end{lemma}
\begin{proof}
  By symmetry we may assume that $m(p)=\dist(a_{\leq d},p)$. Define a template $t$ by
  \[
    t_i\eqdef
    \begin{cases}
      p_i&\text{if }i\leq d\text{ and }a_i\neq p_i,\\
      \star&\text{otherwise}.
    \end{cases}
  \]
  The template $t$ is regular with respect to $a$, of weight $\dist(a_{\leq d},p)$. Since
  $S^{(\ell)}(p)\subset S^{(\ell)}[t]$, it follows from \Cref{Sk-template-bound} that
  \[
    \abs{S^{(\ell)}(p)}\leq \sum_{k\leq d} \abs{S_k^{(\ell)}[t]}\leq \sum_{k\leq d} Z(n/2)^{k-m(p)-L}\binom{2d-m(p)-L}{k-m(p)-L}
                                            \leq d Z (n/2)^{d-m(p)}\binom{2d-m(p)}{d-m(p)}.\qedhere
  \]
\end{proof}
Let $D_m\eqdef \{p\in [n]^d : m(p)\leq m\}$. Note that
\begin{align*}
  \abs{D_m}&\leq \binom{d}{d-m}\binom{m}{d-m}n^{2m-d}
\end{align*}
since every word $p\in D_m$ must agree both with $a_{\leq d}$ and $b_{\leq d}$ in at least $d-m$ positions each.

From the preceding lemma it follows that
\begin{align*}
  \abs{S^{(\ell)}}&=\sum_{p\in [n]^d} \abs{S^{(\ell)}(p)}\leq \sum_m \abs{D_m\setminus D_{m-1}}  \cdot dZ(n/2)^{d-m}\binom{2d-m}{d-m}\\
         &\leq \sum_m \binom{d}{d-m}\binom{m}{d-m}n^{2m-d}\cdot dZ(n/2)^{d-m}\binom{2d-m}{d-m}.\\
\intertext{Making the substitution $t=d-m$ and using the inequality  $\binom{x}{t}\leq \left(\frac{ex}{t}\right)^d$, we obtain}
  \abs{S}&\leq dZn^d\sum_t \binom{d}{t}\binom{d-t}{t}\binom{d+t}{t}(2n)^{-t}\\
         &\leq dZn^d\sum_t \left(\frac{e^3 d(d-t)(d+t)}{t^3}\right) (2n)^{-t}\\
         &\leq d^2Zn^d\max_t\left(\frac{e^3d^3}{2nt^3}\right)^t.\\\intertext{Since the maximum is attained at $t=d(2n)^{-1/3}$, it follows that}
  \abs{S}&\leq d^2Zn^de^{3d(2n)^{-1/3}}\leq d^2Zn^d\left(1+8n^{-1/3}\right)^d.
\end{align*}
In view of the bound for $Z=\Sun(p+n,2L)$ in \Cref{sunflower}, this completes the proof of \Cref{thm:mainL}. In addition, since $Z=1$ for $L=0$, \Cref{thm:main} also follows.

\section{Small \texorpdfstring{$d$-maximal}{d-maximal} sets: proof of \texorpdfstring{\Cref{thm:lower}}{Theorem 6}}
\paragraph{Preliminaries.} Let $S$ be a $d$-maximal set, and let $m\eqdef \abs{S}$ be its size.
We wish to prove a lower bound on $m$. Clearly we may assume that $S$ is finite, and so we may embed $S$ into $[n]^r$ for some finite~$r$.
The dimension parameter $r$ satisfies $r\geq d$, for otherwise $[n]^r$ contains no subsets of diameter $d$ at all.

\paragraph{Fractional solution.} 
For $i\in [n]$, let $e_i\in \R^n$ be the $i$'th standard basis vector.
To each $a\in [n]^r$ we then associate the $r$-by-$n$ matrix $M(a)$ whose $i$'th row is $M(a)_i\eqdef e_{a_i}$, for all $i\in [r]$.
We call matrices of the form $M(a)$ \emph{permissible}. Permissible matrices are precisely
the $\{0,1\}$-matrices of dimension $r$-by-$n$ that have exactly one $1$ in every row.

The Hamming distance on $[n]^r$ and the $\ell_1$ distance on $\R^{r\times n}$ are related by
\[
  \norm{M(a)-M(b)}_1=2\dist(a,b)\qquad\text{ for all }a,b\in [n]^r.
\]

Put $A^{(0)}\eqdef \frac{1}{m}\sum_{a\in S} M(a)$ and observe that, for $b\in S$, we have
\begin{align}
\notag  \norm{A^{(0)}-M(b)}_1&=\frac{1}{m}\,\,\,\sum_{a\in S}\,\,\norm{M(a)-M(b)}_1 \\
\notag                      &=\frac{1}{m}\sum_{a\in S\setminus\{b\}}\norm{M(a)-M(b)}_1\\
\notag                      &=\frac{2}{m}\sum_{a\in S\setminus\{b\}}\dist(a,b)\\
\label{fract-dist}          &\leq \frac{2}{m}\cdot (m-1)d,
\end{align}
where the last inequality follows from the assumption $\diam(S)\leq d$. Hence, we can think of $A^{(0)}$ as a fractional solution
to the problem of finding an element that is close to all points of $S$.

\paragraph{Preparing for the rounding.}
Leaning on the assumption that $m$ is small, we shall round the matrix $A^{(0)}$ to a permissible matrix $A$ that is
close in the $\ell_1$ norm to every $M(b)$, but that is not in the $M$-image of $S$. This will contradict maximality of~$S$.

Consider the following constraints on the matrix $A^{(1)}\in \R^{Mn}$:
\begin{subequations}
\begin{align}
 \label{cond-row-sums}  &\sum_{j=1}^n A^{(1)}_{i,j}=1&&\text{for all }i\in [r],\\
 \label{cond-dist}  &\norm{A^{(1)}-M(b)}_1=\norm{A^{(0)}-M(b)}_1&&\text{for all }b\in S,\\
 \label{cond-cube}  &0\leq A^{(1)}_{i,j}\leq 1&&\text{for all }i\in[r],j\in[n].
\end{align}
\end{subequations}
The conditions in \eqref{cond-row-sums} are linear in the entries of~$A^{(1)}$. The conditions in \eqref{cond-dist} are not linear for general matrices,
but they are linear if we restrict to matrices with entries in the interval $[0,1]$. We therefore obtain a system of $r+m$ linear equations in entries of $A^{(1)}$.

\paragraph{Floating coloring argument.} We treat entries of $A^{(1)}$ as variables. Initially all variables are ``floating'', and their
values are same as the corresponding entries of $A^{(0)}$. Once a variable becomes equal to $0$ or $1$, it is ``frozen'', and its value will not ever change again.
The floating variables thus have values in $(0,1)$. There are $m$ constraints of the form \eqref{cond-dist}, and initially there are
$r$ constraints of the form \eqref{cond-row-sums}. However, as soon as all the variables in a row are frozen, the corresponding constraint is dropped.

Let $f$ be the number of rows that contain at least one floating variable. As long as the number of floating variables is less than $f+m$,
the solution space to linear system made by \eqref{cond-row-sums} and \eqref{cond-dist}
contains a line, which intersects the cube $[0,1]^{r\times n}$ in a boundary point. This corresponds to a solution in which one more of the floating variables
is frozen. We can then update $A^{(1)}$ to such a solution and continue the process.

The process stops when the number of floating variables becomes at most $f+m$. Since row sums are equal to $1$, each row containing a floating variable
in fact contains at least two floating variables. Hence $2f\leq f+m$.

\paragraph{Random rounding.} Since each row of $A^{(1)}$ sums to $1$, we interpret them as probabilities,
and we define $c_i$ to be $j$ with probability $B_{i,j}$. The resulting vector $c=(c_1,\dotsc,c_r)$ is a random
element of~$[n]^r$.

Fix $b\in S$. Since the number of coordinates $c_i$ in $c$ that are not deterministic is equal to $f\leq m$ and $\norm{M(c)-M(b)}_1$ is a $2$-Lipschitz function
of $c$, Azuma's inequality (see \cite[Theorem~7.2.1]{alon_spencer}) implies that
\[
  \Pr\Bigl[\norm{M(c)-M(b)}_1-\norm{A^{(1)}-M(b)}> 2\lambda \sqrt{m}\Bigr]< \exp(-\lambda^2/2).
\]
Substituting \eqref{fract-dist} and \eqref{cond-dist} into this inequality and using the union bound we obtain
\[
  \Pr\Bigl[\dist(c,b)>\frac{m-1}{m}d+\lambda\sqrt{m}\text{ for some }b\in S\Bigr]<m\exp(-\lambda^2/2).
\]

Set $\lambda=\sqrt{2\log m}$. If $d/m- \sqrt{2m\log m}\geq 1$, then there is a choice of $c$ satisfying $\dist(c,b)\leq d-1$ for all $b\in S$.
By maximality of $S$ and the triangle inequality, it follows that the $1$-ball around $c$ is contained in $S$. Since $r\geq d$, it then follows that
$\abs{S}\geq d$. If $d/m- \sqrt{2m\log m}\leq 1$, then $m\geq \frac{1}{2}(d/\log 2d)^{2/3}$.\smallskip

\textsc{Remark.} For $n=2$, the bound in \Cref{thm:lower} can be improved to $m\geq cd^{2/3}$ by replacing, in the proof above,
the random rounding step by Spencer's entropy argument \cite{spencer_six}, or one of its subsequent refinements. For example,
the improvement follows by iteratively applying Theorem 5.1 in \cite{bansal_constructive}. It is likely that this extends to 
$n>2$ as well. For example, Spencer's entropy argument was worked out for a slightly different notion of multicolor discrepancy in \cite{doerr_srivastav}.
However, we believe that the exponent of $2/3$ is not sharp, and so opted for simpler exposition over proving a slightly stronger result.

\section{\texorpdfstring{$d$-maximal}{d-maximal} sets over binary alphabet: proofs of \texorpdfstring{\Cref{thm:binaryconstr,thm:binary}}{Theorems 2 and 3}}\label{sec:binary}
\paragraph{Better construction (proof of \texorpdfstring{\Cref{thm:binaryconstr}}{Theorem 2}).} We exhibit a $d$-maximal set of size~$\binom{\lfloor 3d/2\rfloor}{d}$ in~$\twoinfty$.
In the construction we identify elements of $\twoinfty$ with subsets of $\N$ in the natural way.\smallskip

Suppose that $d=2k$ is even. Let $\binom{[3k]}{2k}\eqdef \{A\subset [3k]:|A|=2k\}$ and set $S_{\text{even}}\eqdef \{\emptyset\}\cup \binom{[3k]}{2k}$.
It is clear that $S_{\text{even}}$ is of diameter $2k$, and we claim that $S_{\text{even}}$ is in fact $2k$-maximal.
Indeed, assume that $S\cup\{B\}$ is also of diameter $2k$. By permuting the elements of $[3k]$ we may assume
that $B\cap [3k]=[m]$ for some integer $m$. We cannot have $m>2k$, for otherwise $\dist(\emptyset,B)\geq m>2k$.
We cannot have $0<m<2k$ either for otherwise
\begin{align*}
  \dist([3k]\setminus [k],B)&\geq \dist([3k]\setminus [k],B\cap [2k])=\min(4k-m,2k+m)>2k.
\end{align*}
So, $B\in S_{\text{even}}$.\smallskip

Suppose that $d=2k+1$ is odd. Let $S_{\text{odd}}\eqdef \{A\subset [3k+1] : A\cap [3k]\in S_{\text{even}}\}$. The inequality $\diam(S_{\text{odd}})\leq 2k+1$
follows from $\diam(S_{\text{even}})\leq 2k$. As to maximality, assume that $S\cup\{B\}$ is also of diameter $2k+1$, and $B\cap [3k]=[m]$.
Let $X\eqdef \{3k+1\}\setminus B$. As above, we cannot have
$m>2k$, for otherwise $\dist(X,B)>2k+1$. We cannot have $0<m<2k$ either, for otherwise $\dist([3k]\cup X\setminus [k],B)=1+\min(4k-m,2k+m)>2k+1$.
So, $B\in S_{\text{odd}}$.

\paragraph{Better upper bound (proof of \texorpdfstring{\Cref{thm:binary}}{Theorem 3}).} 
We continue using the notation from \Cref{sec:basic}; namely, an element $w\in S$ is fixed arbitrarily, and $S_k\eqdef \{a\in S : \dist(a,w)=k\}$.
The alphabet size is $n=2$ throughout this section. For $i\in \N$, let
\[
  p_i\eqdef \Pr_{a\sim S_k}[a_i\neq w_i],
\]
where the randomness is over a choice of a word $a$ uniformly at random from~$S_k$.

\begin{lemma}
  Let $p_1,p_2,p_3,\dotsc$ be as above. Then
  \begin{align}
    \sum_i p_i&=k.\label{eq:psum}\\
    \sum_i p_i^2&\geq k-d/2.\label{eq:psq}
  \end{align}
\end{lemma}  
\begin{proof}
  Choose $a,b$ independently uniformly at random from~$S_k$. Because
  $\sum_i p_i=\E_{a\sim S_k} \dist(a,w)$, the definition of $S_k$ implies \eqref{eq:psum}.
  Similarly, since $S_k$ is of diameter at most $d$, we deduce that
  \[
    d\geq \E_{a,b\sim S_k}[\dist(a,b)]=\sum_i 2p_i(1-p_i).
  \]
  Taken together with \eqref{eq:psum}, this implies \eqref{eq:psq}.
\end{proof}

\begin{lemma}\label{above-half}
  Let $\ell$ be a positive integer, and suppose that $\Delta \eqdef \sum_{i\leq \ell} (p_i-\tfrac{1}{2})$ satisfies $\Delta\geq 0$. Then
  \[
    \abs{S_k}\leq  \binom{2d}{k}(1+\Delta/\ell)^{-\Delta}.
  \]
\end{lemma}
\begin{proof}
  Set $m\eqdef \lceil \Delta\rceil$. 
  For a word $a\in \Ninfty$, let $a_{\leq \ell}\in [n]^{\ell}$ be its prefix of length~$\ell$.
  Pick a set $X$ of $m$ elements in $[\ell]$ uniformly among all $\binom{\ell}{m}$ such sets.
  Let $t$ be the template with $t_i=1$ for $i\in X$ and $t_i=\star$ for $i\notin X$.
  Note that, for $a\in \Ninfty$,
  \[
    \Pr[a\text{ fits }t]=
    \binom{\dist(a_{\leq \ell},w_{\leq \ell})}{m}\Big/\binom{\ell}{m}.
  \]
  The function
  \begin{equation*}
    f(x)\eqdef
    \begin{cases}
      \binom{x}{m}/\binom{\ell}{m}&\text{if }x\geq m-1,\\
      0&\text{if }x\leq m-1
    \end{cases}
  \end{equation*}
  is convex on all of $\R$. Therefore, 
  \begin{align*}
    \E\bigl[\abs{S_k[t]}\bigr]&=\sum_{a\in S_k}\Pr[a\text{ fits }t]=\sum_{a\in S_k} f\bigl(\dist(a_{\leq \ell},w_{\leq \ell})\bigr)\mkern-88mu\\
                              &\geq \abs{S_k}f\left(\E_{a\sim S_k} \dist(a_{\leq \ell},w_{\leq \ell})\right)\\
                              &=\abs{S_k}f(\Delta+\ell/2)\\
                              &=\abs{S_k}\binom{\Delta+\ell/2}{m}\Big/\binom{\ell}{m}&&\text{because }\Delta+\ell/2\geq m\\
                              &= 2^{-m}\abs{S_k}\prod_{i\in [m]}\frac{2(\Delta+\ell/2-m+i)}{\ell-m+i}\\
                              &= 2^{-m}\abs{S_k}\prod_{i\in [m]}\left(1+\frac{\Delta+i-1}{\ell-m+i}\right)&&\text{because }m\leq \Delta+1\\
                              &\geq 2^{-m}\abs{S_k}\prod_{i\in [m]}\left(1+\frac{\Delta}{\ell}\right).
  \end{align*}
  Since $t$ is a regular template of weight $m$, it follows from \Cref{Sk-template-bound} that $\abs{S_k[t]}\leq \binom{2d-m}{k-m}\leq \left(\frac{k}{2d}\right)^m\binom{2d}{k}\leq 2^{-m}\binom{2d}{k}$, and
  so
  \begin{equation*}
    \abs{S_k}\leq \binom{2d}{k}(1+\Delta/\ell)^{-m}.
  \end{equation*}
  Since $m\geq \Delta$, the result follows.
\end{proof}

Denote by $\Ball_n(r)$ a Hamming ball of radius $r$ around some point of $[2]^n$.
\begin{lemma}[Kleitman's theorem \cite{kleitman}] 
  Suppose that $d<n$. Then
  \begin{enumerate}
  \item
    For even $d$, the set $\Ball_n(d/2)$ is a largest diameter-$d$ set in $[2]^n$.
  \item
    For odd $d$, the set $[2]\times \Ball_{n-1}\bigl((d-1)/2\bigr)$ is a largest diameter-$d$ set in $[2]^n$.
  \end{enumerate}
\end{lemma}
In fact, Kleitman proved that these examples are unique, though we will not rely on this.

\begin{corollary}
  Assume that $d< n$. Then every diameter-$d$ subset of $[2]^n$ is of size at most $2^{n H(d/2n)+1}$, where $H(\cdot)$ is the binary entropy function.
\end{corollary}
\begin{proof}
  By \cite[Exercise 15.8.5]{alon_spencer}, the size of $\Ball_n(r)$ is bounded by $2^{nH(r/n)}$ for $r\leq n/2$.
\end{proof}

Let $S_k^{\leq \ell}\eqdef \{a_{\leq \ell} : a\in S_k\}$. The preceding corollary implies that
\[
  \abs{S_k^{\leq \ell}}\leq 2^{\ell H(d/2\ell)+1}.
\]

We can assemble these ingredients together to prove \Cref{thm:binary}. By permuting the coordinates, we may assume that $p_1\geq p_2\geq p_3\geq \dotsb$.
Let $\ell\eqdef 2d$, and define $\Delta$ as in \Cref{above-half}.
Let $\Pbig\eqdef \sum_{i\leq 2d} p_i$ and $\Psmall\eqdef \sum_{i>2d}$. Observe that
\begin{align*}
  \sum_{i\leq 2d} p_i^2&=\sum_{i\in [\ell]} (p_i-1/2)p_i+\Pbig/2\leq \Delta+\Pbig/2,\\
  \sum_{i>2d} p_i^2&\leq p_{2d+1}\Psmall\leq \frac{\Pbig}{2d}\Psmall.
\end{align*}
Adding these up and using \eqref{eq:psum} and \eqref{eq:psq} yields
\begin{align*}
  k-d/2&\leq \Delta+\Pbig/2+\frac{\Pbig}{2d}\Psmall\\
       &\leq \Delta+\Pbig/2+\frac{d-\Psmall}{2d}\Psmall\\
       &=\Delta+k/2-\Psmall^2/2d,\\\intertext{and so}
  \Psmall^2&\leq d^2-dk+2d\Delta. 
\end{align*}
So, either $\Delta>d/40000$, in which case we are done by \Cref{above-half}, or
\[
  \Psmall\leq d/100 \qquad\text{ for }k\geq 1999d/20000.
\]

For a word $s\in \Ninfty$, let $s_{>\ell}$ be its suffix obtained by deleting the first $\ell$ symbols.
Since $\Psmall=\E_{a\sim S_k}[\dist(a_{>\ell},w_{>\ell})]$, Markov's inequality implies that 
the set $S_k'\eqdef \{a\in S : \dist(a_{>\ell},w_{>\ell}) \leq d/50\}$ is of size
\[
  \abs{S_k'}\geq \abs{S_k}/2.
\]

Let $X\eqdef \{a_{\leq 2d} : a\in S_k'\}$. Since $X$ is of diameter at most $d$,
and so $\abs{X}\leq 2^{2d H(1/4)+1}\leq 2^{1.63d+1}$. By \Cref{Sk-template-bound} it then follows that, for $k\geq 19999d/20000$,
\[
  \abs{S_k'}\leq \sum_{x\in X} \abs{\{a\in S_k' : a_{\leq 2d}=x\}}\leq \sum_{x\in X}\binom{2d-\dist(x,w_{\leq 2d})}{d/50}\leq 2^{1.63d+1}\binom{2d}{d/50}\leq 2^{1.83d+1}.
\]
Hence,
\[
  \abs{S}\leq \sum_{k<199/200d} \abs{S_k}+\sum_{k\geq 19999/20000d}\abs{S_k}\leq \sum_{k\leq 19999/20000d}\binom{2d}{k}+d\cdot 2^{1.83d+2}\leq (4-10^{-10})^d
\]
for all large~$d$.

\bibliographystyle{plain}
\bibliography{refs.bib}

\end{document}